\documentclass{amsart}

\usepackage{fancyhdr}
\usepackage{lastpage}
\usepackage{yhmath}

\newcommand{\bdis}{\begin{displaymath}}
\newcommand{\edis}{\end{displaymath}}
\newcommand{\be}{\begin{equation}}
\newcommand{\ee}{\end{equation}}

\newcommand{\mcal}{\mathcal}

\newcommand{\vp}{\varphi}
\newcommand{\vt}{\vartheta}

\newcommand{\zf}{\zeta\left(\frac{1}{2}+it\right)}

\DeclareMathOperator*{\ssum}{\sum\sum}

\theoremstyle{definition}

\theoremstyle{remark}
\newtheorem{remark}[]{Remark}

\newtheorem*{mydef1}{{\bf Theorem}}

\newtheorem*{mydef51}{{\bf Lemma 1}}

\newtheorem*{mydef52}{{\bf Lemma 2}}

\newtheorem*{mydef53}{{\bf Lemma 3}}

\newtheorem*{mydef5A}{{\bf Lemma A}}
\newtheorem*{mydef5B}{{\bf Lemma B}}

\newtheorem*{mydef6}{\bf Formula}

\numberwithin{equation}{section}

\begin{document}

\title{$\Omega$-theorem for short trigonometric sum}

\author{Jan Moser}

\address{Department of Mathematical Analysis and Numerical Mathematics, Comenius University, Mlynska Dolina M105, 842 48 Bratislava, SLOVAKIA}

\email{jan.mozer@fmph.uniba.sk}


\begin{abstract}
We obtain in this paper new application of the classical E.C. Titchmarsh' discrete method (1934) in the theory of the Riemann
$\zf$ - function. Namely, we shall prove the first localized $\Omega$-theorem for short trigonometric sum. This paper is the English version
of the work of reference \cite{4}.
\end{abstract}

\maketitle

\section{Result}

\subsection{}

In this paper we shall study the following short trigonometric sum
\be
S(t,T,K)=\sum_{e^{-1/K}P_0<n<P_0}\cos(t\ln n),\
P_0=\sqrt{\frac{T}{2\pi}},\ t\in [T,T+U],
\ee
where
\be
U=T^{1/2}\psi\ln T,\ \psi\leq K\leq T^{1/6}\ln^2T,\ \psi <\ln T,
\ee
and $\psi=\psi(T)$ stands for arbitrary slowly increasing function unbounded (from above). For example
\bdis
\psi(T)=\ln\ln T,\ \ln\ln\ln T,\ \dots
\edis
Let
\bdis
\{ t_\nu\}
\edis
be the Gram-Titchmarsh sequence defined by the formula
\bdis
\vt(t_\nu)=\pi\nu,\quad \nu=1,2,\dots
\edis
(see \cite{6}, pp. 221, 329) where
\bdis
\vt(t)=-\frac t2\ln\pi+\mbox{Im}\left\{\ln\Gamma\left(\frac 14+i\frac t2\right)\right\}.
\edis
Next, we denote by
\bdis
G(T,K,\psi)
\edis
the number of such $t_\nu$ that (see (1.1)) obey the following
\be
t_\nu\in [T,T+U] \ \wedge \ |S(t_\nu,T,K)|>\frac 12 \sqrt{\frac{P_0}{K}}=AT^{1/4}K^{-1/2}.
\ee
The following theorem holds true.

\begin{mydef1}
There are
\bdis
T_0(K,\psi)>0, \ A>0
\edis
such that
\be
G(T,K,\psi)>AT^{1/6}K^{-1}\psi \ln^2T,\ T\geq T_0(K,\psi).
\ee
\end{mydef1}

\subsection{}

Let us remind the following estimate by Karatsuba (see \cite{1}, p. 89)
\bdis
\overset{*}{S}(x)=\sum_{1\leq n\leq x}n^{it}=\mcal{O}(\sqrt{x}t^{\epsilon}),\ 0<x<t
\edis
holds true on the Lindel\" of hypothesis ($0<\epsilon$ is an arbitrary small number in this). Of course,
\be
\begin{split}
& \overset{*}{S}(t,T,K)=\sum_{e^{-1/K}P_0\leq n\leq P_0}n^{it}=\mcal{O}(\sqrt{P_0}T^\epsilon)=\mcal{O}(T^{1/4+\epsilon}), \\
& t\in [T,T+U],
\end{split}
\ee
and
\bdis
|\overset{*}{S}(t,T,K)|\geq |S(t,T,K)|.
\edis

\begin{remark}
 Since (see (1.1), (1.3), (1.4)) the inequality
 \bdis
 |S(t,T,\ln\ln T)|> A\frac{T^{1/4}}{\sqrt{\ln\ln T}},\ T\to\infty
 \edis
 is fulfilled for arbitrary big $t$, then the Karatsuba's estimate (1.5) is an almost exact estimate.
\end{remark}

\subsection{}

With regard to connection between short trigonometric sum and the theory of the Riemann zeta function see our paper \cite{3}.

\section{Main lemmas and proof of Theorem}

Let
\be
w(t,T,K)=\sum_{e^{-1/K}P_0<n<P_0}\frac{1}{\sqrt{n}}\cos(t\ln n),
\ee
\be
w_1(t,T,K)=\sum_{e^{-1/K}P_0<n<P_0}\left(\frac{1}{\sqrt{n}}-\frac{1}{\sqrt{P_0}}\right)\cos(t\ln n).
\ee
The following lemmas hold true.

\begin{mydef5A}
 \be
 \sum_{T\leq t_\nu\leq T+U}w^2(t_\nu,T,K)=\frac{1}{4\pi}UK^{-1}\ln\frac{T}{2\pi}+\mcal{O}(K^{-1}\sqrt{T}\ln^2T).
 \ee
\end{mydef5A}

\begin{mydef5B}
 \be
 \sum_{T\leq t_\nu\leq T+U}w_1^2(t_\nu,T,K)=\frac{1}{48\pi}UK^{-3}\ln\frac{T}{2\pi}+\mcal{O}(K^{-3}\sqrt{T}\ln^2T).
 \ee
\end{mydef5B}

\begin{remark}
 Of course, the formulae (2.3), (2.4) are asymptotic ones (see (1.2)).
\end{remark}

Now we use the main lemmas A and B for completion of the Theorem. Since (see \cite{2}, (23))
\be
Q_0=\sum_{T\leq t_\nu\leq T+U}1=\frac{1}{2\pi}U\ln\frac{T}{2\pi}+\mcal{O}\frac{U^2}{T},
\ee
then we obtain (see (2.3), (2.4)) that
\be
\frac{1}{Q_0}\sum_{T\leq t_\nu\leq T+U}w^2(t_\nu,T,K)\sim \frac{1}{2K},\ T\to\infty,
\ee
\be
\frac{1}{Q_0}\sum_{T\leq t_\nu\leq T+U}w_1^2(t_\nu,T,K)\sim \frac{1}{24K^3},
\ee
\be
\frac{1}{Q_0}\sum_{T\leq t_\nu\leq T+U}w\cdot w_1=\mcal{O}\left(\frac{1}{K^2}\right),
\ee
(we used the Schwarz inequality in (2.8)). Next, we have (see (1.1), (2.1), (2.2)) that
\be
w(t_\nu,T,K)=\frac{1}{\sqrt{P_0}}S(t_\nu,T,K)+w_1(t_\nu,T,K).
\ee
Consequently we obtain (see (2.6) -- (2.9)) the following

\begin{mydef6}
 \be
 \frac{1}{Q_0}\sum_{T\leq t_\nu\leq T+U}S^2\sim \frac{P_0}{2K},\ T\to\infty.
 \ee
\end{mydef6}
Next, we denote by $Q_1$ the number of such values
\bdis
t_\nu\in [T,T+U],
\edis
that fulfill the inequality (see (1.3))
\be
|S|>\frac{1}{2}\sqrt{\frac{P_0}{K}};\quad Q_1=G(T,K,\psi),
\ee
and
\bdis
Q_0-Q_1=Q_2
\edis
(see (2.5)). Since (see (1.1) and \cite{6}, p. 92)
\bdis
|S(t,T,K)|<A \sqrt{P_0}T^{1/6},\ t\in [T,T+U],
\edis
then we have (see (2.10), (2.11)) that
\bdis
\frac{1}{3K}< AT^{1/3}\frac{Q_1}{Q_0}+\frac{1}{4K}\frac{Q_2}{Q_0}<A T^{1/3}\frac{Q_1}{Q_0}+\frac{1}{4K},
\edis
i. e.
\be
AQ_1>\frac{1}{12}Q_0T^{-1/3}K^{-1}.
\ee
Consequently, we obtain (see (1.2), (2.5), (2.11), (2.12)) the following estimate
\bdis
Q_1=G>A T^{1/6}K^{-1}\psi \ln^2T
\edis
that is required result (1.4).

\section{Lemma 1}

Let
\be
w_2=\ssum_{e^{-1/K}P_0<n<m<P_0}\frac{1}{\sqrt{nm}}\cos\left( t_\nu\ln\frac nm\right).
\ee

The following lemma holds true.

\begin{mydef51}
 \be
 \sum_{T\leq t_\nu\leq T+U}w_2=\mcal{O}(K^{-1}\sqrt{T}\ln^2T).
 \ee
\end{mydef51}

\begin{proof}
 The following inner sum (comp. \cite{5}, p. 102; $t_{\nu+1}\longrightarrow t_\nu$)
 \be
 w_{21}=\sum_{T\leq T_\nu\leq T+U}\cos\{ 2\pi\psi_1(\nu)\},
 \ee
 where
 \bdis
 \psi_1(\nu)=\frac{1}{2\pi}t_\nu\ln\frac nm
 \edis
 applies to our sum (3.2).
 Now we obtain by method \cite{5}, pp. 102-103 the following estimate
 \be
 w_{21}=\mcal{O}\left(\frac{\ln T}{\ln\frac nm}\right).
 \ee
Since (see (1.2))
\bdis
e^{-1/K}>1-\frac 1K\geq 1-\frac{1}{\psi}>\frac 12,
\edis
then
\bdis
2m>2e^{-1/K}P_0>P_0,\quad m\in (e^{-1/K}P_0,P_0),
\edis
i. e. in our case (see (3.1)) we have that
\bdis
2m>n.
\edis
Consequently, the method \cite{6}, p. 116, $\sigma=\frac 12,\ m=n-r$ gives the estimate
\be
\begin{split}
 & \ssum_{e^{-1/K}P_0<n<m<P_0}\frac{1}{\sqrt{mn}\ln\frac nm}< \\
 & <
 A\sum_{e^{-1/K}P_0<n<P_0}\sum_{r\leq n/2}\frac 1r< AK^{-1}P_0\ln P_0<AK^{-1}\sqrt{T}\ln T,
\end{split}
\ee
where
\be
\sum_{e^{-1/K}P_0<n<P_0} 1\sim \frac{P_0}{K}.
\ee
Now, required result (3.2) follows from (3.1), (3.3) -- (3.5).
\end{proof}

\section{Lemma 2}

Let
\be
w_3=\ssum_{e^{-1/K}P_0<m<n<P_0}\frac{1}{\sqrt{mn}}\cos\{ t_\nu\ln(mn)\}.
\ee
The following lemma holds true.
\begin{mydef52}
 \be
 \sum_{T\leq t_\nu\leq T+U}w_3=\mcal{O}(K^{-1}\sqrt{T}\ln^2T).
 \ee
\end{mydef52}

\begin{proof}
 The following inner sum (comp. \cite{5}, p. 103; $t_{\nu+1}\longrightarrow t_\nu$)
 \bdis
 w_{31}=\sum_{T\leq t_\nu\leq T+U}\cos\{ 2\pi \chi(\nu)\},
 \edis
 where
 \bdis
 \chi(\nu)=\frac{1}{2\pi}t_\nu\ln(nm)
 \edis
 applies to our sum (4.2).  Next, the method \cite{5}, pp. 103-104 gives us that
 \be \begin{split}
 & w_{31}=\int_{\chi'(x)<1/2}\cos\{ 2\pi\chi(x)\}{\rm d}x+ \\
  & +  \int_{\chi'(x)>1/2}\cos[2\pi\{ \chi(x)-x\}]{\rm d}x+\mcal{O}(1)=J_1+J_2+\mcal{O}(1),
\end{split} \ee
where
\bdis
J_1=\mcal{O}\left(\frac{\ln T}{\ln n}\right)=\mcal{O}(1),\ n\in (e^{-1/K}P_0,P_0),
\edis
and ($m<n<2m,\ n=m+r$)
\bdis
J_2=\mcal{O}\left( \frac{m\ln(m+1)}{r}\right).
\edis
Now, the term $J_1$ contributes to the sum (4.2), (comp. (3.6)) as
\be
\begin{split}
 & \mcal{O}\left(\ssum_{e^{-1/K}P_0<m<n<P_0}\frac{1}{\sqrt{mn}}\right) = \\
 & = \mcal{O}\left( \frac{1}{P_0}\ssum_{e^{-1/K}P_0<m<n<P_0} 1\right)= \\
 & = \mcal{O}\left(\frac{1}{P_0}\frac{P_0^2}{K^2}\right)=\mcal{O}(K^{-2}\sqrt{T}),
\end{split}
\ee
and the same contribution corresponds to the term $\mcal{O}(1)$ in (4.3), while the contribution of the term $J_2$ is
\be
\begin{split}
 & \mcal{O}\left(\sum_{e^{-1/K}P_0<m<P_0}\frac{1}{\sqrt{m}}\sum_{r=1}^m\frac{1}{\sqrt{m}}\frac{m\ln(m+1)}{r}\right)= \\
 & = \mcal{O}\left(\frac{P_0}{K}\ln^2P_0\right)=\mcal{O}(K^{-1}\sqrt{T}\ln^2T).
\end{split}
\ee
Now, the required result (4.2) follows from (4.1), (4.4), (4.5).
\end{proof}

\section{Lemma 3}

Let
\be
w_4=\sum_{e^{-1/K}P_0<n<P_0}\frac 1n.
\ee
The following lemma holds true.

\begin{mydef53}
\be
\sum_{T\leq t_\nu\leq T+U}w_4=\mcal{O}(K^{-1}\sqrt{T}\ln^2T).
\ee
\end{mydef53}

\begin{proof}
The following inner sum
\bdis
w_{41}=\sum_{T\leq t_\nu\leq T+U}\cos\{ 2\pi\chi_1(\nu)\},
\edis
where
\bdis
\chi_1(\nu)=\frac{1}{\pi}t_\nu\ln n.
\edis
applies to our sum (5.2).  Since (comp. \cite{5}, p. 103)
\bdis
\chi_1'(\nu)=\frac{\ln n}{\vt'(t_\nu)},
\edis
next, (comp. \cite{5}, p. 100)
\bdis
\begin{split}
 & \vt'(t_\nu)=\frac 12\ln\frac{t_\nu}{2\pi}+\mcal{O}\left(\frac{1}{t_\nu}\right)=
 \frac 12\ln\frac{T}{2\pi}+\mcal{O}\left(\frac UT\right)+\mcal{O}\left(\frac 1T\right)\sim \\
 & \sim \ln P_0,
\end{split}
\edis
and
\bdis
\ln P_0-\frac 1K<\ln n<\ln P_0,\ n\in (e^{-1/K}P_0,P_0),
\edis
then
\bdis
\chi_1'(\nu)\sim 1.
\edis
Since, for example,
\bdis
\frac 12<\chi_1'(\nu)<\frac 32,
\edis
then we have (comp. \cite{5}, p. 104) that
\bdis
w_{41}=\int\cos[2\pi\{\chi_1(x)-x\}]{\rm d}x+\mcal{O}(1)=J_3+\mcal{O}(1).
\edis
Now, (comp. \cite{5}, p. 104)
\bdis
\chi_1''(\nu)<-A\frac{\ln n}{T\ln^3T}<-\frac{B}{T\ln^2T},\ n\in (e^{-1/K}P_0,P_0),
\edis
and
\bdis
J_3=\mcal{O}(\sqrt{T}\ln T),\quad w_{41}=\mcal{O}(\sqrt{T}\ln T).
\edis
Consequently, we get the required result (5.2)
\bdis
\sum_{T\leq t_\nu\leq T+U}w_4=\mcal{O}\left(\sqrt{T}\ln T\sum_{e^{-1/K}P_0<n<P_0}\frac 1n\right)=
\mcal{O}(K^{-1}\sqrt{T}\ln T),
\edis
where
\bdis
\sum_{e^{-1/K}P_0<n<P_0}\frac 1n\sim \frac 1K
\edis
by the well-known Euler's formula
\bdis
\sum_{1\leq n<x}\frac 1n=\ln x+c+\mcal{O}\left(\frac 1x\right),
\edis
where $c$ is the Euler's constant.
\end{proof}

\section{Lemmas A and B}

\subsection{Proof of Lemma A}

First of all, we have (see (2.1)) that
\be
\begin{split}
& w^2(t_\nu,T,K)=\\
& = \ssum_{e^{-1/K}P_0<m,n<P_0}\frac{1}{\sqrt{mn}}\cos( t_\nu\ln m)\cos( t_\nu\ln n)= \\
& = \frac 12 \sum_n \frac 1n+\ssum_{m<n}\frac{1}{\sqrt{mn}}\cos\left( t_\nu\ln\frac nm\right)+ \\
& + \ssum_{m<n}\frac{1}{\sqrt{mn}}\cos\{ t_\nu\ln (mn)\}+\frac 12\sum_{n}\frac 1n\cos( 2t_\nu\ln n)= \\
& = \frac{1}{2K}+\mcal{O}\left(\frac{1}{\sqrt{T}}\right)+w_2+w_3+w_4,
\end{split}
\ee
(see (5.3), (3.1), (4.1), (5.1)). Consequently, we obtain the required result (2.3) from (6.1) by (2.5), (3.2), (4.2), (5.2).

\subsection{Proof of Lemma B} First of all we have (see (2.2)) that
\bdis
w_1(t_\nu,T,K)=\sum_{e^{-1/K}P_0<n<P_0}\frac{\alpha(n)}{\sqrt{n}}\cos(t_\nu\ln n),
\edis
where
\bdis
\alpha(n)=1-\sqrt{\frac{n}{P_0}}.
\edis
Of course, $\alpha(n)$ is decreasing and
\be
0<\alpha(n)<\frac 1K,\quad  n\in (e^{-1/K}P_0,P_0).
\ee
Next, (comp. (6.1))
\be
\begin{split}
& w_1^2(t_\nu,T,K)= \\
& = \frac 12\sum_n\frac{\alpha^2(n)}{n}+\ssum_{m<n}\frac{\alpha(m)\alpha(n)}{\sqrt{mn}}\cos\left( t_\nu\ln\frac nm\right)+ \\
& + \ssum_{m<n}\frac{\alpha(m)\alpha(n)}{\sqrt{mn}}\cos\left( t_\nu\ln(mn)\right)+\frac 12\sum_n
\frac{\alpha^2(n)}{n}\cos( 2t_\nu\ln n)= \\
& = \frac 12\bar{w}_1+\bar{w}_2+\bar{w}_3+\frac 12\bar{w}_4.
\end{split}
\ee
Since (see (6.2))
\bdis
\alpha(m)\alpha(n)<K^{-2}
\edis
then we obtain by a similar way as in the case of the estimates (3.2), (4.2) and (5.2) that
\be
\sum_{T\leq t_\nu\leq T+U}\left\{ \bar{w}_2+\bar{w}_3+\frac 12\bar{w}_4\right\}=
\mcal{O}(K^{-3}\sqrt{T}\ln^2T).
\ee
In the case of the sum
\be
\frac 12\sum_{T\leq t_\nu\leq T+U}\bar{w}_1
\ee
we use the following summation formula (see \cite{6}, p. 13)
\bdis
\begin{split}
& \sum_{a\leq n<b}\vp(n)=\int_a^b \vp(x){\rm d}x+\int_a^b
\left( x-[x]-\frac 12\right)\vp'(x){\rm d}x+ \\
& + \left( a-[a]-\frac 12\right)\vp(a)-\left( b-[b]-\frac 12\right)\vp(b)
\end{split}
\edis
in the case
\bdis
a=e^{-1/K}P_0, b=P_0, \vp(x)=\frac{\alpha^2(x)}{x}=\frac 1x-\frac{2}{\sqrt{P_0x}}+\frac{1}{P_0}.
\edis
Hence
\bdis
\begin{split}
& \int_{e^{-1/K}P_0}^{P_0}\frac{\alpha^2(x)}{x}{\rm d}x=\frac{1}{K}-4\left( 1-e^{-1/(2K)}\right)+1-e^{-1/K}= \\
& = \frac{1}{12K^3}+\mcal{O}(K^{-4}),
\end{split}
\edis
and
\bdis
\vp'(x)=\mcal{O}(P_0^{-2}),\
\vp(e^{-1/K}P_0)=\mcal{O}\left(\frac{1}{x^2P_0}\right),\ \vp(P_0)=0.
\edis
Consequently, we have (see (6.5))
\bdis
\frac 12\bar{w}_1=\frac{1}{24K^3}+\mcal{O}\left(\frac{1}{K^4}\right),
\edis
and (see (1.2), (2.5))
\be
\begin{split}
& \frac 12\sum_{T\leq t_\nu\leq T+U}\bar{w}_1=\frac{1}{48\pi}UK^{-3}\ln\frac{T}{2\pi}+
\mcal{O}\left(\frac{U\ln T}{K^4}\right)+\mcal{O}\left(\frac{U^2}{K^3T}\right)= \\
& = \frac{1}{48\pi}UK^{-3}\ln\frac{T}{2\pi}+\mcal{O}(K^{-3}\sqrt{T}\ln T).
\end{split}
\ee
Finally, we obtain the required result (2.4) from (6.3) by (6.4) -- (6.6).

\end{document}